\documentclass{amsart}
\usepackage{eurosym}
\usepackage{amsmath}
\usepackage{amsfonts}

\newtheorem{theorem}{Theorem}
\theoremstyle{plain}

\newtheorem{definition}{Definition}

\newtheorem{lemma}{Lemma}

\newtheorem{remark}{Remark}

\numberwithin{equation}{section}
\input{tcilatex}

\begin{document}
\title[On Hermite-Hadamard-Fej\'{e}r type Inequality for Convex...]{On
Hermite-Hadamard-Fej\'{e}r type Inequality for Convex Functions via
Fractional Integrals}
\author{Abdullah AKKURT}
\address{[Department of Mathematics, Faculty of Science and Arts, University
of Kahramanmara\c{s} S\"{u}t\c{c}\"{u} \.{I}mam, 46100, Kahramanmara\c{s},
Turkey}
\email{abdullahmat@gmail.com}
\author{H\"{u}seyin YILDIRIM}
\address{[Department of Mathematics, Faculty of Science and Arts, University
of Kahramanmara\c{s} S\"{u}t\c{c}\"{u} \.{I}mam, 46100, Kahramanmara\c{s},
Turkey}
\email{hyildir@ksu.edu.tr}
\keywords{\textbf{}Integral inequalities, Fractional integrals,
Hermite-Hadamard-Fej\'{e}r Inequality.}

\begin{abstract}
In this paper, we have established some generalized integral inequalities of
Hermite-Hadamard-Fej\'{e}r type for generalized fractional integrals. The
results presented here would provide generalizations of those given in
earlier works.
\end{abstract}

\maketitle

\section{Introduction}

Let $f:I\subset 
\mathbb{R}
\rightarrow 
\mathbb{R}
$ be a convex function define on an interval $I$ of real numbers, and $%
a,b\in I$ with $a<b.$ Then the following inequalities hold:

\begin{equation}
\begin{array}{c}
f\left( \dfrac{a+b}{2}\right) \leq \dfrac{1}{b-a}\int\limits_{a}^{b}f\left(
x\right) dx\leq \dfrac{f\left( a\right) +f\left( b\right) }{2}.%
\end{array}
\label{1.1}
\end{equation}

It was first discovered by Hermite in $1881$ in the Journal Mathesis. This
inequality (\ref{1.1}) was nowhere mentioned in the mathematical literature
untill $1893$. In \cite{4}, Beckenbach, a leading expert on the theory of
convex functions, wrote that the inequality (\ref{1.1}) was proved by
Hadamard in $1893$. In $1974$, Mitrinovi\v{c} found Hermite and Hadamard's
note in Mathesis. That is why, the inequality (\ref{1.1}) was known as
Hermite-Hadamard inequality. We note that Hermite-Hadamard's inequality may
be regarded as a refinements of the concept of convexity and it follows
easily from Jensen's inequality. This inequality (\ref{1.1}) has been
received renewed attention in recent years and a remarkable variety of
refinements and generalizations have been found in \cite{4}-\cite{16}, \cite%
{18}, \cite{20}.

The most well known inequalities connected with the integral mean of a
convex functions are Hermite-Hadamard inequalities or its weighted versions,
the so-called Hermite-Hadamard-Fej\v{e}r inequality.

In \cite{10}, Fej\'{e}r established the following Fej\'{e}r inequality which
is the weighted generalization of Hermite-Hadamard inequalities (\ref{1.1}):

\begin{theorem}
Let $f:I\rightarrow 
\mathbb{R}
$ be a convex on $I$ and let $a,b\in I$ with $a<b$. Then the inequality%
\begin{equation}
\begin{array}{c}
f\left( \dfrac{a+b}{2}\right) \int\limits_{a}^{b}g(x)dx\leq \dfrac{1}{b-a}%
\int\limits_{a}^{b}f(t)g(x)dt\leq \dfrac{f(a)+f(b)}{2}\int%
\limits_{a}^{b}g(x)dx.%
\end{array}
\label{1.2}
\end{equation}%
holds, where $f:\left[ a,b\right] \rightarrow 
\mathbb{R}
$ is nonnegative, integrable, and symmetric to $\dfrac{a+b}{2}$.
\end{theorem}

In \cite{11}, Sarikaya et al. represented Hermite--Hadamard's inequalities
in fractional integral forms as follows.

\begin{theorem}
Let $f:\left[ a,b\right] \rightarrow 
\mathbb{R}
$ be a positive function with $0\leq a<b$ and $f\in L\left[ a,b\right] .$ If 
$f$ is a convex function on $\left[ a,b\right] ,$ then the following
inequalities for fractional integrals holds%
\begin{equation}
f\left( \frac{a+b}{2}\right) \leq \frac{\Gamma \left( \alpha +1\right) }{%
2(b-a)^{\alpha }}\left[ J_{a^{+}}^{\alpha }f(b)+J_{b^{-}}^{\alpha }f(a)%
\right] \leq \frac{f(a)+f(b)}{2},  \label{1.3}
\end{equation}%
with $\alpha >0.$
\end{theorem}

In \cite{5} Set et. al. obtained the following lemma.

\begin{lemma}
Let $f:\left[ a,b\right] \rightarrow 
\mathbb{R}
$ be a differentiable mapping on $\left( a,b\right) $ with $a<b$ and let $g:%
\left[ a,b\right] \rightarrow 
\mathbb{R}
$. If $f^{\prime }$, $g\in L\left[ a,b\right] $, then the following identity
for fractional integrals holds:%
\begin{equation}
\begin{array}{l}
f\left( \frac{a+b}{2}\right) \left[ J_{\left( \frac{a+b}{2}\right)
^{-}}^{\alpha }g(a)+J_{\left( \frac{a+b}{2}\right) ^{+}}^{\alpha }g(b)\right]
\\ 
\ \ \ \ -\left[ J_{\left( \frac{a+b}{2}\right) ^{-}}^{\alpha }\left(
gf\right) (a)+J_{\left( \frac{a+b}{2}\right) ^{+}}^{\alpha }\left( fg\right)
(b)\right] =\dfrac{1}{\Gamma \left( \alpha \right) }\int%
\limits_{a}^{b}k(t)df(t),%
\end{array}
\label{1.4}
\end{equation}%
where%
\begin{equation*}
\begin{array}{c}
k(t)=\left\{ 
\begin{array}{cc}
\int\limits_{a}^{t}(s-a)^{\alpha -1}g(s)ds & t\in \left[ a,\frac{a+b}{2}%
\right] , \\ 
\int\limits_{b}^{t}(b-s)^{\alpha -1}g(s)ds & t\in \left[ \frac{a+b}{2},b%
\right] .%
\end{array}%
\right.%
\end{array}%
\end{equation*}
\end{lemma}

We give some necessary definitions and mathematical preliminaries of
fractional calculus theory which are used throughout this paper.

\begin{definition}
Let $h\left( x\right) $ be an increasing and positive monotone function on $%
[0,\infty )$, also derivative $h^{^{\prime }}\left( x\right) $ is continuous
on $[0,\infty )\ $and $h\left( 0\right) =0$. The space $X_{h}^{p}\left(
0,\infty \right) \left( 1\leq p<\infty \right) $ of those real-valued
Lebesque measurable functions $f$ on $[0,\infty )\ $for which%
\begin{equation}
\begin{array}{cc}
\left\vert \left\vert f\right\vert \right\vert _{X_{h}^{p}} & =\left(
\int\limits_{0}^{\infty }\left\vert f(t)\right\vert ^{p}h^{^{\prime }}\left(
x\right) dt\right) ^{\frac{1}{p}}<\infty ,\ \ \ \ 1\leq p\leq \infty%
\end{array}
\label{1.5}
\end{equation}%
and for the case $p=\infty $%
\begin{equation}
\begin{array}{cc}
\left\vert \left\vert f\right\vert \right\vert _{X_{h}^{\infty }} & 
=ess\sup\limits_{1\leq t<\infty }\left[ f(t)h^{^{\prime }}\left( x\right) %
\right] .%
\end{array}
\label{1.6}
\end{equation}
\end{definition}

\begin{definition}
(\cite{6}). In particular, when $h\left( x\right) =x\ \left( 1\leq p<\infty
\right) $ the space $X_{h}^{p}\left( 0,\infty \right) \ $coincides with the $%
L_{p}[0,\infty )-$space ($\left\vert \left\vert f\right\vert \right\vert
_{X_{h}^{\infty }}=\left\vert \left\vert f\right\vert \right\vert _{\infty }$%
) and also if we take $h\left( x\right) =\dfrac{x^{k+1}}{k+1}\ (k\geq 0)\ $%
the space $X_{h}^{p}\left( 0,\infty \right) \ $coincides with the $%
L_{p,k}[0,\infty )-$space.
\end{definition}

\begin{definition}
(\cite{1}).$\ $Let $\left( a,b\right) $ be a finite interval of the real
line $%
\mathbb{R}
$ and $\alpha >0$. Also let $h\left( x\right) $ be an increasing and
positive monotone function on $(a,b]$, having a continuous derivative $%
h^{^{\prime }}\left( x\right) $ on $(a,b)$. The left- and right-sided
fractional integrals of a function $\ f$ \ with respect to another function $%
h$ on $[a,b]$ are defined by%
\begin{equation}
\begin{array}{c}
\left( J_{a^{+},h}^{\alpha }f\right) \left( x\right) :=\frac{1}{\Gamma
\left( \alpha \right) }\int_{a}^{x}\left[ h\left( x\right) -h\left( t\right) %
\right] ^{\alpha -1}h^{^{\prime }}\left( t\right) f\left( t\right) dt,\
x\geq a%
\end{array}
\label{1.7}
\end{equation}%
and%
\begin{equation}
\begin{array}{c}
\left( J_{b^{-},h}^{\alpha }f\right) \left( x\right) :=\frac{1}{\Gamma
\left( \alpha \right) }\int_{x}^{b}\left[ h\left( t\right) -h\left( x\right) %
\right] ^{\alpha -1}h^{^{\prime }}\left( t\right) f\left( t\right) dt,\
x\leq b.%
\end{array}
\label{1.8}
\end{equation}
\end{definition}

\begin{definition}
If we take $h(x)=x$, then the equalities (\ref{1.7}) and (\ref{1.8}) will be

\begin{equation}
\begin{array}{c}
\left( J_{a^{+}}^{\alpha }f\right) (x)=\frac{1}{\Gamma (\alpha )}%
\int\limits_{a}^{x}(x-t)^{\alpha -1}f(t)dt,\ x>a%
\end{array}
\label{1.9}
\end{equation}%
and%
\begin{equation}
\begin{array}{c}
\left( J_{b^{-}}^{\alpha }f\right) (x)=\frac{1}{\Gamma (\alpha )}%
\int\limits_{x}^{b}(t-x)^{\alpha -1}f(t)dt,\ b>x.%
\end{array}
\label{1.10}
\end{equation}%
These integrals are called left-sided\ Riemann-Liouville fractional integral
and right-sided\ Riemann-Liouville fractional integral respectively \cite{1}-%
\cite{3}, \cite{6}, \cite{17}.
\end{definition}

In this paper, we have established some generalized fractional integral
inequalities. The results presented here would provide generalizations of
those given in earlier works.

\section{Main Results}

\begin{lemma}
\label{lem1} Let $f:\left[ a,b\right] \rightarrow 
\mathbb{R}
$ be a differentiable mapping on $\left( a,b\right) $ with $a<b$ and let $g:%
\left[ a,b\right] \rightarrow 
\mathbb{R}
$. If $f^{\prime }$, $g\in X_{h}^{p}\left[ a,b\right] $, then the following
identity for fractional integrals holds:%
\begin{equation}
\begin{array}{l}
f\left( h\left( \frac{a+b}{2}\right) \right) \left[ J_{\left( \frac{a+b}{2}%
\right) ^{-}}^{\alpha }g(h(a))+J_{\left( \frac{a+b}{2}\right) ^{+}}^{\alpha
}g(h(b))\right] \\ 
\\ 
\ \ \ \ -\left[ J_{\left( \frac{a+b}{2}\right) ^{-}}^{\alpha }\left( g\times
\left( f\circ h\right) \right) (a)+J_{\left( \frac{a+b}{2}\right)
^{+}}^{\alpha }\left( g\times \left( f\circ h\right) \right) (b)\right] \\ 
\\ 
=\dfrac{1}{\Gamma \left( \alpha \right) }\int\limits_{a}^{b}k(t)df(h(t))%
\end{array}
\label{2.1}
\end{equation}%
where%
\begin{equation*}
\begin{array}{c}
k(t)=\left\{ 
\begin{array}{cc}
\int\limits_{a}^{t}(h(s)-h(a))^{\alpha -1}g(s)h^{\prime }(s)ds & t\in \left[
a,\frac{a+b}{2}\right] , \\ 
&  \\ 
\int\limits_{b}^{t}(h(b)-h(s))^{\alpha -1}g(s)h^{\prime }(s)ds & t\in \left[ 
\frac{a+b}{2},b\right] .%
\end{array}%
\right.%
\end{array}%
\end{equation*}

\begin{proof}
It suffices to note that%
\begin{equation*}
\begin{array}{l}
I=\int\limits_{a}^{b}k(t)df(h(t)) \\ 
\\ 
=\int\limits_{a}^{\frac{a+b}{2}}\left(
\int\limits_{a}^{t}(h(s)-h(a))^{\alpha -1}g(s)h^{\prime }(s)ds\right)
df(h(t)) \\ 
\\ 
\ \ \ \ +\int\limits_{\frac{a+b}{2}}^{b}\left(
\int\limits_{b}^{t}(h(b)-h(s))^{\alpha -1}g(s)h^{\prime }(s)ds\right)
df(h(t)) \\ 
\\ 
=I_{1}+I_{2}.%
\end{array}%
\end{equation*}%
By integration by parts, we get%
\begin{equation*}
\begin{array}{l}
I_{1}=\int\limits_{a}^{\frac{a+b}{2}}\left(
\int\limits_{a}^{t}(h(s)-h(a))^{\alpha -1}g(s)h^{\prime }(s)ds\right)
df(h(t)) \\ 
=\left. \left( \int\limits_{a}^{t}(h(s)-h(a))^{\alpha -1}g(s)h^{\prime
}(s)ds\right) f(h(t))\right\vert _{a}^{\frac{a+b}{2}} \\ 
\ \ \ \ -\int\limits_{a}^{\frac{a+b}{2}}(h(t)-h(a))^{\alpha
-1}g(t)f(h(t))h^{\prime }(t)dt%
\end{array}%
\end{equation*}%
\begin{equation*}
\begin{array}{l}
=f\left( h\left( \frac{a+b}{2}\right) \right) \int\limits_{a}^{\frac{a+b}{2}%
}(h(s)-h(a))^{\alpha -1}g(s)h^{\prime }(s)ds \\ 
\ \ \ \ -\int\limits_{a}^{\frac{a+b}{2}}(h(t)-h(a))^{\alpha
-1}g(t)f(h(t))h^{\prime }(t)dt \\ 
=\Gamma \left( \alpha \right) \left[ f\left( h\left( \frac{a+b}{2}\right)
\right) J_{\left( \frac{a+b}{2}\right) ^{-}}^{\alpha }g(h(a))-J_{\left( 
\frac{a+b}{2}\right) ^{-}}^{\alpha }\left( g\left( f\circ h\right) \right)
(a)\right] 
\end{array}%
\end{equation*}%
and similarly%
\begin{equation*}
\begin{array}{l}
I_{2}=\int\limits_{\frac{a+b}{2}}^{b}\left(
\int\limits_{b}^{t}(h(b)-h(s))^{\alpha -1}g(s)h^{\prime }(s)ds\right)
df(h(t)) \\ 
=\left. \left( \int\limits_{b}^{t}(h(s)-h(a))^{\alpha -1}g(s)h^{\prime
}(s)ds\right) f(h(t))\right\vert _{\frac{a+b}{2}}^{b} \\ 
\ \ \ \ -\int\limits_{\frac{a+b}{2}}^{b}(h(b)-h(t))^{\alpha
-1}g(t)f(h(t))h^{\prime }(t)dt \\ 
=f\left( h\left( \frac{a+b}{2}\right) \right) \int\limits_{\frac{a+b}{2}%
}^{b}(h(b)-h(s))^{\alpha -1}g(s)h^{\prime }(s)ds \\ 
\ \ \ \ -\int\limits_{\frac{a+b}{2}}^{b}(h(b)-h(t))^{\alpha
-1}g(t)f(h(t))h^{\prime }(t)dt \\ 
=\Gamma \left( \alpha \right) \left[ f\left( h\left( \frac{a+b}{2}\right)
\right) J_{\left( \frac{a+b}{2}\right) ^{+}}^{\alpha }g(h(b))-J_{\left( 
\frac{a+b}{2}\right) ^{+}}^{\alpha }\left( g\left( f\circ h\right) \right)
(b)\right] .%
\end{array}%
\end{equation*}%
Thus, can write 
\begin{equation*}
\begin{array}{l}
I=I_{1}+I_{2} \\ 
=\Gamma \left( \alpha \right) \left\{ f\left( h\left( \frac{a+b}{2}\right)
\right) \left[ J_{\left( \frac{a+b}{2}\right) ^{-}}^{\alpha
}g(h(a))+J_{\left( \frac{a+b}{2}\right) ^{+}}^{\alpha }g(h(b))\right]
\right.  \\ 
\ \ \ \ \left. -\left[ J_{\left( \frac{a+b}{2}\right) ^{-}}^{\alpha }\left(
g\left( f\circ h\right) \right) (a)+J_{\left( \frac{a+b}{2}\right)
^{+}}^{\alpha }\left( g\left( f\circ h\right) \right) (b)\right] \right\} .%
\end{array}%
\end{equation*}%
Multiplying the both sides $\left( \Gamma \left( \alpha \right) \right)
^{-1},\ $we obtain (\ref{2.1}) which complates the proof.
\end{proof}
\end{lemma}

\begin{remark}
If we choose $h(x)=x$ in Lemma \ref{lem1}, then the inequality (\ref{2.1})
reduces to $\left( 1.4\right) .$
\end{remark}

\begin{remark}
If we choose $h(x)=x$,\ $g(x)=1\ $and $\alpha =1\ $in Lemma \ref{lem1}, we
obtain Lemma 2.1 in \cite{21}.

\begin{theorem}
\label{thm1} Let $f:I\rightarrow 
\mathbb{R}
$ be a differentiable mapping on $I^{\circ }$ and $f^{\prime }\in X_{h}^{p}%
\left[ a,b\right] $ with $a<b$ and $g:\left[ a,b\right] \rightarrow 
\mathbb{R}
$ is continuous. If $\left\vert f^{\prime }\right\vert \ $is convex on $%
\left[ a,b\right] ,\ $then the following inequality for fractional integrals
holds:%
\begin{equation}
\begin{array}{l}
\left\vert f\left( h\left( \frac{a+b}{2}\right) \right) \left[ J_{\left( 
\frac{a+b}{2}\right) ^{-}}^{\alpha }g(h(a))+J_{\left( \frac{a+b}{2}\right)
^{+}}^{\alpha }g(h(b))\right] \right.  \\ 
\\ 
\ \ \ \ -\left. \left[ J_{\left( \frac{a+b}{2}\right) ^{-}}^{\alpha }\left(
g\times \left( f\circ h\right) \right) (a)+J_{\left( \frac{a+b}{2}\right)
^{+}}^{\alpha }\left( g\times \left( f\circ h\right) \right) (b)\right]
\right\vert  \\ 
\\ 
\leq \frac{\left\vert \left\vert g\right\vert \right\vert _{X_{h}^{\infty }%
\left[ a,\frac{a+b}{2}\right] ,\infty }}{\left( h(b)-h(a)\right) \Gamma
\left( \alpha +1\right) }\left\{ \left\vert f^{\prime }\left( h(a)\right)
\right\vert \left[ \frac{(h(\frac{a+b}{2})-h(a))^{\alpha +1}}{\alpha +1}%
\left( h(b)-h(a)\right) \right. \right.  \\ 
\\ 
\ \ \ \ -\left. \left\vert f^{\prime }\left( h(a)\right) \right\vert \frac{%
(h(\frac{a+b}{2})-h(a))^{\alpha +2}}{\alpha +2}\right] 
\end{array}
\label{2.2}
\end{equation}%
\begin{equation*}
\begin{array}{l}
\ \ \ \ \left. +\left\vert f^{\prime }\left( h(b)\right) \right\vert \left[ 
\frac{(h(\frac{a+b}{2})-h(a))^{\alpha +2}}{\alpha +2}\right] \right\} +\frac{%
\left\vert \left\vert g\right\vert \right\vert _{X_{h}^{\infty }\left[ \frac{%
a+b}{2},b\right] ,\infty }}{\left( h(b)-h(a)\right) \Gamma \left( \alpha
+1\right) } \\ 
\\ 
\times \left\{ \left\vert f^{\prime }\left( h(b)\right) \right\vert \left[ 
\frac{(h(b)-h(\frac{a+b}{2}))^{\alpha +1}}{\alpha +1}\left( h(b)-h(a)\right)
-\frac{(h(b)-h(\frac{a+b}{2}))^{\alpha +2}}{\alpha +2}\right] \right.  \\ 
\\ 
\ \ \ \ \left. +\left\vert f^{\prime }\left( h(a)\right) \right\vert \left[ 
\frac{(h(b)-h(\frac{a+b}{2}))^{\alpha +2}}{\alpha +2}\right] \right\} 
\end{array}%
\end{equation*}%
with $\alpha >0.$

\begin{proof}
If $\left\vert f^{\prime }\right\vert \ $is convex on $\left[ a,b\right] ,\ $%
we know that for $t\in \left[ a,b\right] $%
\begin{equation*}
\begin{array}{l}
\left\vert f^{\prime }(h(t))\right\vert  \\ 
\\ 
=\left\vert f^{\prime }\left( \dfrac{h(b)-h(t)}{h(b)-h(a)}h(a)+\dfrac{%
h(t)-h(a)}{h(b)-h(a)}h(b)\right) \right\vert  \\ 
\\ 
\leq \dfrac{h(b)-h(t)}{h(b)-h(a)}\left\vert f^{\prime }\left( h(a)\right)
\right\vert +\dfrac{h(t)-h(a)}{h(b)-h(a)}\left\vert f^{\prime }\left(
h(b)\right) \right\vert .%
\end{array}%
\end{equation*}%
From Lemma \ref{lem1}, we have%
\begin{equation*}
\begin{array}{l}
\left\vert f\left( h\left( \frac{a+b}{2}\right) \right) \left[ J_{\left( 
\frac{a+b}{2}\right) ^{-}}^{\alpha }g(h(a))+J_{\left( \frac{a+b}{2}\right)
^{+}}^{\alpha }g(h(b))\right] \right.  \\ 
\\ 
\ \ \ \ \left. -\left[ J_{\left( \frac{a+b}{2}\right) ^{-}}^{\alpha }\left(
g\left( f\circ h\right) \right) (a)+J_{\left( \frac{a+b}{2}\right)
^{+}}^{\alpha }\left( g\left( f\circ h\right) \right) (b)\right] \right\vert 
\end{array}%
\end{equation*}%
\begin{equation*}
\begin{array}{l}
\leq \frac{1}{\Gamma \left( \alpha \right) }\left\{ \int\limits_{a}^{\frac{%
a+b}{2}}\left\vert \int\limits_{a}^{t}(h(s)-h(a))^{\alpha -1}g(s)h^{\prime
}(s)ds\right\vert \left\vert f^{\prime }\left( h(t)\right) \right\vert
h^{\prime }(t)dt\right.  \\ 
\ \ \ \ \left. +\int\limits_{\frac{a+b}{2}}^{b}\left\vert
\int\limits_{b}^{t}(h(b)-h(s))^{\alpha -1}g(s)h^{\prime }(s)ds\right\vert
\left\vert f^{\prime }\left( h(t)\right) \right\vert h^{\prime
}(t)dt\right\}  \\ 
\leq \frac{\left\vert \left\vert g\right\vert \right\vert _{X_{h}^{\infty }%
\left[ a,\frac{a+b}{2}\right] ,\infty }}{\left( h(b)-h(a)\right) \Gamma
\left( \alpha \right) }\left\{ \int\limits_{a}^{\frac{a+b}{2}}\left\vert
\int\limits_{a}^{t}(h(s)-h(a))^{\alpha -1}h^{\prime }(s)ds\right\vert
\right.  \\ 
\ \ \ \ \times \left. \left( h(b)-h(t)\left\vert f^{\prime }\left(
h(a)\right) \right\vert \right) h^{\prime }(t)dt\right. 
\end{array}%
\end{equation*}%
\begin{equation*}
\begin{array}{l}
\ \ \ \ \left. +\int\limits_{a}^{\frac{a+b}{2}}\left\vert
\int\limits_{a}^{t}(h(s)-h(a))^{\alpha -1}h^{\prime }(s)ds\right\vert \left(
h(t)-h(a)\left\vert f^{\prime }\left( h(b)\right) \right\vert \right)
h^{\prime }(t)dt\right\}  \\ 
\ \ \ \ +\frac{\left\vert \left\vert g\right\vert \right\vert
_{X_{h}^{\infty }\left[ \frac{a+b}{2},b\right] ,\infty }}{\left(
h(b)-h(a)\right) \Gamma \left( \alpha \right) }\left\{ \int\limits_{\frac{a+b%
}{2}}^{b}\left\vert \int\limits_{b}^{t}(h(b)-h(s))^{\alpha -1}h^{\prime
}(s)ds\right\vert \right.  \\ 
\ \ \ \ \times \left. \left( h(b)-h(t)\left\vert f^{\prime }\left(
h(a)\right) \right\vert \right) h^{\prime }(t)dt\right.  \\ 
\ \ \ \ \left. +\int\limits_{\frac{a+b}{2}}^{b}\left\vert
\int\limits_{b}^{t}(h(b)-h(s))^{\alpha -1}h^{\prime }(s)ds\right\vert \left(
h(t)-h(a)\left\vert f^{\prime }\left( h(b)\right) \right\vert \right)
h^{\prime }(t)dt\right\} .%
\end{array}%
\end{equation*}%
\begin{equation*}
\begin{array}{l}
\leq \frac{\left\vert \left\vert g\right\vert \right\vert _{X_{h}^{\infty }%
\left[ a,\frac{a+b}{2}\right] ,\infty }}{\left( h(b)-h(a)\right) \Gamma
\left( \alpha +1\right) }\left\{ \int\limits_{a}^{\frac{a+b}{2}%
}(h(t)-h(a))^{\alpha }\left( h(b)-h(t)\left\vert f^{\prime }\left(
h(a)\right) \right\vert \right) h^{\prime }(t)dt\right.  \\ 
\left. +\int\limits_{a}^{\frac{a+b}{2}}(h(t)-h(a))^{\alpha +1}\left\vert
f^{\prime }\left( h(b)\right) \right\vert h^{\prime }(t)dt\right\}  \\ 
+\frac{\left\vert \left\vert g\right\vert \right\vert _{X_{h}^{\infty }\left[
\frac{a+b}{2},b\right] ,\infty }}{\left( h(b)-h(a)\right) \Gamma \left(
\alpha +1\right) }\left\{ \int\limits_{\frac{a+b}{2}}^{b}(h(b)-h(t))^{\alpha
+1}\left\vert f^{\prime }\left( h(a)\right) \right\vert h^{\prime
}(t)dt\right.  \\ 
\left. +\int\limits_{\frac{a+b}{2}}^{b}(h(b)-h(t))^{\alpha -1}\left(
h(t)-h(a)\left\vert f^{\prime }\left( h(b)\right) \right\vert \right)
h^{\prime }(t)dt\right\} .%
\end{array}%
\end{equation*}%
\begin{equation*}
\begin{array}{l}
\leq \frac{\left\vert \left\vert g\right\vert \right\vert _{X_{h}^{\infty }%
\left[ a,\frac{a+b}{2}\right] ,\infty }}{\left( h(b)-h(a)\right) \Gamma
\left( \alpha +1\right) }\left\{ \left\vert f^{\prime }\left( h(a)\right)
\right\vert \left[ \frac{(h(\frac{a+b}{2})-h(a))^{\alpha +1}}{\alpha +1}%
\left( h(b)-h(a)\right) \right. \right.  \\ 
\ \ \ \ \left. -\frac{(h(\frac{a+b}{2})-h(a))^{\alpha +2}}{\alpha +2}\right] 
\\ 
\ \ \ \ \left. +\left\vert f^{\prime }\left( h(b)\right) \right\vert \left[ 
\frac{(h(\frac{a+b}{2})-h(a))^{\alpha +2}}{\alpha +2}\right] \right\} 
\end{array}%
\end{equation*}%
\begin{equation*}
\begin{array}{l}
\ \ \ \ +\frac{\left\vert \left\vert g\right\vert \right\vert
_{X_{h}^{\infty }\left[ \frac{a+b}{2},b\right] ,\infty }}{\left(
h(b)-h(a)\right) \Gamma \left( \alpha +1\right) }\left\{ \left\vert
f^{\prime }\left( h(b)\right) \right\vert \left[ \frac{(h(b)-h(\frac{a+b}{2}%
))^{\alpha +1}}{\alpha +1}\left( h(b)-h(a)\right) \right] \right.  \\ 
\ \ \ \ \left. -\frac{(h(b)-h(\frac{a+b}{2}))^{\alpha +2}}{\alpha +2}\right. 
\\ 
\ \ \ \ \left. +\left\vert f^{\prime }\left( h(a)\right) \right\vert \left[ 
\frac{(h(b)-h(\frac{a+b}{2}))^{\alpha +2}}{\alpha +2}\right] \right\} .%
\end{array}%
\end{equation*}%
This completes the proof.
\end{proof}
\end{theorem}
\end{remark}

\begin{remark}
If we choose $h(x)=x\ $in Theorem \ref{thm1}, we obtain Theorem $6$ in \cite%
{5}.
\end{remark}

\begin{remark}
If we choose $h(x)=x$,\ $g(x)=1\ $and $\alpha =1\ $in Theorem \ref{thm1}, we
obtain Theorem $2.2$ in \cite{21}.

\begin{theorem}
\label{thm2} Let $f:I\rightarrow 
\mathbb{R}
$ be a differentiable mapping on $I^{\circ }$ and $f^{\prime }\in X_{h}^{p}%
\left[ a,b\right] $ with $a<b$ and $g:\left[ a,b\right] \rightarrow 
\mathbb{R}
$ is continuous. If $\left\vert f^{\prime }\right\vert ^{q}\ $is convex on $%
\left[ a,b\right] ,\ q\geq 1\ $then the following inequality for fractional
integrals holds:%
\begin{equation}
\begin{array}{l}
f\left( h\left( \frac{a+b}{2}\right) \right) \left[ J_{\left( \frac{a+b}{2}%
\right) ^{-}}^{\alpha }g(h(a))+J_{\left( \frac{a+b}{2}\right) ^{+}}^{\alpha
}g(h(b))\right]  \\ 
\\ 
\ \ \ \ -\left[ J_{\left( \frac{a+b}{2}\right) ^{-}}^{\alpha }\left( g\times
\left( f\circ h\right) \right) (a)+J_{\left( \frac{a+b}{2}\right)
^{+}}^{\alpha }\left( g\times \left( f\circ h\right) \right) (b)\right]  \\ 
\\ 
\leq \frac{\left\vert \left\vert g\right\vert \right\vert _{X_{h}^{\infty }%
\left[ a,\frac{a+b}{2}\right] ,\infty }}{\left[ \left( h(b)-h(a)\right) %
\right] ^{\frac{1}{q}}\Gamma \left( \alpha \right) }\left( \frac{h\left( 
\frac{a+b}{2}\right) -h(a)}{\alpha \left( \alpha +1\right) }\right) ^{1-%
\frac{1}{q}}\left\{ \left\vert f^{\prime }\left( h(a)\right) \right\vert
^{q}\right.  \\ 
\\ 
\ \ \ \ \times \left. \left[ \frac{(h(\frac{a+b}{2})-h(a))^{\alpha +1}\left(
h(b)-h(a)\right) }{\alpha +1}-\frac{(h(\frac{a+b}{2})-h(a))^{\alpha +2}}{%
\alpha +2}\right] \right. 
\end{array}
\label{2.3}
\end{equation}%
\begin{equation*}
\begin{array}{l}
\ \ \ \ \left. +\frac{\left\vert f^{\prime }\left( h(b)\right) \right\vert
^{q}(h(\frac{a+b}{2})-h(a))^{\alpha +2}}{\alpha +2}\right\} ^{\frac{1}{q}}
\\ 
\\ 
\ \ \ \ +\frac{\left\vert \left\vert g\right\vert \right\vert
_{X_{h}^{\infty }\left[ \frac{a+b}{2},b\right] ,\infty }}{\left[ \left(
h(b)-h(a)\right) \right] ^{\frac{1}{q}}\Gamma \left( \alpha \right) }\left( 
\frac{h\left( b\right) -h(\frac{a+b}{2})}{\alpha \left( \alpha +1\right) }%
\right) ^{1-\frac{1}{q}}\left\{ \left\vert f^{\prime }\left( h(b)\right)
\right\vert ^{q}\right.  \\ 
\\ 
\ \ \ \ \times \left. \left[ \frac{(h(b)-h(\frac{a+b}{2}))^{\alpha +1}\left(
h(b)-h(a)\right) }{\alpha +1}-\frac{(h(b)-h(\frac{a+b}{2}))^{\alpha +2}}{%
\alpha +2}\right] \right.  \\ 
\\ 
\ \ \ \ \left. +\left\vert f^{\prime }\left( h(a)\right) \right\vert ^{q}
\left[ \frac{(h(b)-h(\frac{a+b}{2}))^{\alpha +2}}{\alpha +2}\right] \right\}
^{\frac{1}{q}}.%
\end{array}%
\end{equation*}%
with $\alpha >0.$

\begin{proof}
If $\left\vert f^{\prime }\right\vert ^{q}\ $is convex on $\left[ a,b\right]
,\ $we know that for $t\in \left[ a,b\right] $%
\begin{equation*}
\begin{array}{l}
\left\vert f^{\prime }(h(t))\right\vert ^{q} \\ 
\\ 
=\left\vert f^{\prime }\left( \frac{h(b)-h(t)}{h(b)-h(a)}h(a)+\frac{h(t)-h(a)%
}{h(b)-h(a)}h(b)\right) \right\vert ^{q} \\ 
\\ 
\leq \frac{h(b)-h(t)}{h(b)-h(a)}\left\vert f^{\prime }\left( h(a)\right)
\right\vert ^{q}+\frac{h(t)-h(a)}{h(b)-h(a)}\left\vert f^{\prime }\left(
h(b)\right) \right\vert ^{q}.%
\end{array}%
\end{equation*}%
From Lemma \ref{lem1}, power mean inequality and the convexity of $%
\left\vert f^{\prime }\right\vert ^{q},\ $it follows that%
\begin{equation*}
\begin{array}{l}
f\left( h\left( \frac{a+b}{2}\right) \right) \left[ J_{\left( \frac{a+b}{2}%
\right) ^{-}}^{\alpha }g(h(a))+J_{\left( \frac{a+b}{2}\right) ^{+}}^{\alpha
}g(h(b))\right]  \\ 
-\left[ J_{\left( \frac{a+b}{2}\right) ^{-}}^{\alpha }\left( g\left( f\circ
h\right) \right) (a)+J_{\left( \frac{a+b}{2}\right) ^{+}}^{\alpha }\left(
g\left( f\circ h\right) \right) (b)\right]  \\ 
\leq \frac{1}{\Gamma \left( \alpha \right) }\left\{ \left( \int\limits_{a}^{%
\frac{a+b}{2}}\left\vert \int\limits_{a}^{t}(h(s)-h(a))^{\alpha
-1}g(s)h^{\prime }(s)ds\right\vert h^{\prime }(t)dt\right) ^{1-1/q}\right. 
\\ 
\left. \times \left( \int\limits_{a}^{\frac{a+b}{2}}\left\vert
\int\limits_{a}^{t}(h(s)-h(a))^{\alpha -1}g(s)h^{\prime }(s)ds\right\vert
\left\vert f^{\prime }\left( h(t)\right) \right\vert ^{q}h^{\prime
}(t)dt\right) ^{1/q}\right\}  \\ 
+\frac{1}{\Gamma \left( \alpha \right) }\left\{ \left( \int\limits_{\frac{a+b%
}{2}}^{b}\left\vert \int\limits_{b}^{t}(h(b)-h(s))^{\alpha -1}g(s)h^{\prime
}(s)ds\right\vert h^{\prime }(t)dt\right) ^{1-1/q}\right.  \\ 
\times \left. \left( \int\limits_{\frac{a+b}{2}}^{b}\left\vert
\int\limits_{b}^{t}(h(b)-h(s))^{\alpha -1}g(s)h^{\prime }(s)ds\right\vert
\left\vert f^{\prime }\left( h(t)\right) \right\vert ^{q}h^{\prime
}(t)dt\right) ^{1/q}\right\} 
\end{array}%
\end{equation*}%
\begin{equation*}
\begin{array}{l}
\leq \frac{\left\vert \left\vert g\right\vert \right\vert _{X_{h}^{\infty }%
\left[ a,\frac{a+b}{2}\right] ,\infty }}{\Gamma \left( \alpha \right) }%
\left\{ \left( \int\limits_{a}^{\frac{a+b}{2}}\left\vert
\int\limits_{a}^{t}(h(s)-h(a))^{\alpha -1}h^{\prime }(s)ds\right\vert
h^{\prime }(t)dt\right) ^{1-1/q}\right.  \\ 
\ \ \ \ \left. \times \left( \int\limits_{a}^{\frac{a+b}{2}}\left\vert
\int\limits_{a}^{t}(h(s)-h(a))^{\alpha -1}h^{\prime }(s)ds\right\vert
\left\vert f^{\prime }\left( h(t)\right) \right\vert ^{q}h^{\prime
}(t)dt\right) ^{1/q}\right\}  \\ 
\ \ \ \ +\frac{\left\vert \left\vert g\right\vert \right\vert
_{X_{h}^{\infty }\left[ \frac{a+b}{2},b\right] ,\infty }}{\Gamma \left(
\alpha \right) }\left\{ \left( \int\limits_{\frac{a+b}{2}}^{b}\left\vert
\int\limits_{b}^{t}(h(b)-h(s))^{\alpha -1}h^{\prime }(s)ds\right\vert
h^{\prime }(t)dt\right) ^{1-1/q}\right. 
\end{array}%
\end{equation*}%
\begin{equation*}
\begin{array}{l}
\ \ \ \ \left. \times \left( \int\limits_{\frac{a+b}{2}}^{b}\left\vert
\int\limits_{b}^{t}(h(b)-h(s))^{\alpha -1}h^{\prime }(s)ds\right\vert
\left\vert f^{\prime }\left( h(t)\right) \right\vert ^{q}h^{\prime
}(t)dt\right) ^{1/q}\right\} . \\ 
\leq \frac{\left\vert \left\vert g\right\vert \right\vert _{X_{h}^{\infty }%
\left[ a,\frac{a+b}{2}\right] ,\infty }}{\left[ \left( h(b)-h(a)\right) %
\right] ^{\frac{1}{q}}\Gamma \left( \alpha \right) }\left\{ \int\limits_{a}^{%
\frac{a+b}{2}}\left\vert \int\limits_{a}^{t}(h(s)-h(a))^{\alpha -1}h^{\prime
}(s)ds\right\vert \right.  \\ 
\ \ \ \ \times \left. \left( h(b)-h(t)\left\vert f^{\prime }\left(
h(a)\right) \right\vert ^{q}\right) h^{\prime }(t)dt\right.  \\ 
\ \ \ \ \left. +\int\limits_{a}^{\frac{a+b}{2}}\left\vert
\int\limits_{a}^{t}(h(s)-h(a))^{\alpha -1}h^{\prime }(s)ds\right\vert \left(
h(t)-h(a)\left\vert f^{\prime }\left( h(b)\right) \right\vert ^{q}\right)
h^{\prime }(t)dt\right\} ^{\frac{1}{q}}%
\end{array}%
\end{equation*}%
\begin{equation*}
\begin{array}{l}
\ \ \ \ +\frac{\left\vert \left\vert g\right\vert \right\vert
_{X_{h}^{\infty }\left[ \frac{a+b}{2},b\right] ,\infty }}{\left[ \left(
h(b)-h(a)\right) \right] ^{\frac{1}{q}}\Gamma \left( \alpha \right) }\left\{
\int\limits_{\frac{a+b}{2}}^{b}\left\vert
\int\limits_{b}^{t}(h(b)-h(s))^{\alpha -1}h^{\prime }(s)ds\right\vert
\right.  \\ 
\ \ \ \ \times \left. \left( h(b)-h(t)\left\vert f^{\prime }\left(
h(a)\right) \right\vert ^{q}\right) h^{\prime }(t)dt\right.  \\ 
\ \ \ \ \left. +\int\limits_{\frac{a+b}{2}}^{b}\left\vert
\int\limits_{b}^{t}(h(b)-h(s))^{\alpha -1}h^{\prime }(s)ds\right\vert \left(
h(t)-h(a)\left\vert f^{\prime }\left( h(b)\right) \right\vert ^{q}\right)
h^{\prime }(t)dt\right\} ^{\frac{1}{q}} \\ 
\leq \frac{\left\vert \left\vert g\right\vert \right\vert _{X_{h}^{\infty }%
\left[ a,\frac{a+b}{2}\right] ,\infty }}{\left[ \left( h(b)-h(a)\right) %
\right] ^{\frac{1}{q}}\Gamma \left( \alpha \right) }\left( \frac{h\left( 
\frac{a+b}{2}\right) -h(a)}{\alpha \left( \alpha +1\right) }\right) ^{1-1/q}
\\ 
\ \ \ \ \times \left\{ \left[ \frac{\left\vert f^{\prime }\left( h(a)\right)
\right\vert ^{q}(h(\frac{a+b}{2})-h(a))^{\alpha +1}\left( h(b)-h(a)\right) }{%
\alpha +1}\right. \right. 
\end{array}%
\end{equation*}%
\begin{equation*}
\begin{array}{l}
\ \ \ \ \left. \left. -\frac{\left\vert f^{\prime }\left( h(a)\right)
\right\vert ^{q}(h(\frac{a+b}{2})-h(a))^{\alpha +2}}{\alpha +2}\right]
+\left\vert f^{\prime }\left( h(b)\right) \right\vert ^{q}\left[ \frac{(h(%
\frac{a+b}{2})-h(a))^{\alpha +2}}{\alpha +2}\right] \right\} ^{\frac{1}{q}}
\\ 
\ \ \ \ +\frac{\left\vert \left\vert g\right\vert \right\vert
_{X_{h}^{\infty }\left[ \frac{a+b}{2},b\right] ,\infty }}{\left[ \left(
h(b)-h(a)\right) \right] ^{\frac{1}{q}}\Gamma \left( \alpha \right) }\left( 
\frac{h\left( b\right) -h(\frac{a+b}{2})}{\alpha \left( \alpha +1\right) }%
\right) ^{1-1/q} \\ 
\ \ \ \ \times \left\{ \left\vert f^{\prime }\left( h(b)\right) \right\vert
^{q}\left[ \frac{(h(b)-h(\frac{a+b}{2}))^{\alpha +1}\left( h(b)-h(a)\right) 
}{\alpha +1}\right. \right.  \\ 
\ \ \ \ \left. -\frac{(h(b)-h(\frac{a+b}{2}))^{\alpha +2}}{\alpha +2}\right]
\left. +\left\vert f^{\prime }\left( h(a)\right) \right\vert ^{q}\left[ 
\frac{(h(b)-h(\frac{a+b}{2}))^{\alpha +2}}{\alpha +2}\right] \right\} ^{%
\frac{1}{q}}.%
\end{array}%
\end{equation*}
\end{proof}
\end{theorem}
\end{remark}

\begin{remark}
If we choose $h(x)=x\ $in Theorem \ref{thm2}, we obtain Theorem $7$ in \cite%
{5}.
\end{remark}

\begin{lemma}
\label{lem2} Let $f:\left[ a,b\right] \rightarrow 
\mathbb{R}
$ be a differentiable mapping on $\left( a,b\right) $ with $a<b$ and let $g:%
\left[ a,b\right] \rightarrow 
\mathbb{R}
$. If $f^{\prime }$, $g\in X_{h}^{p}\left[ a,b\right] $, then the following
identity for fractional integrals holds:%
\begin{equation}
\begin{array}{l}
\left\{ f\left( h\left( \frac{a+b}{2}\right) \right) \left[
-J_{a^{+}}^{\alpha }g(h(b))-J_{b^{-}}^{\alpha }g(h(a))\right. \right.  \\ 
\\ 
\ \ \ \ +\left. J_{\left( \frac{a+b}{2}\right) ^{+}}^{\alpha
}g(h(b))+J_{\left( \frac{a+b}{2}\right) ^{-}}^{\alpha }g(h(a))\right.  \\ 
\\ 
\ \ \ \ +\frac{1}{\Gamma \left( \alpha \right) }\left. \left.
\int\limits_{a}^{b}(h(b)-h(a))^{\alpha -1}g(s)h^{\prime }(s)ds\right]
\right.  \\ 
\\ 
\ \ \ \ +\left. J_{a^{+}}^{\alpha }\left( g\left( f\circ h\right) \right)
(b)+J_{b^{-}}^{\alpha }\left( g\left( f\circ h\right) \right) (a)\right. 
\end{array}
\label{2.4}
\end{equation}%
\begin{equation*}
\begin{array}{l}
\ \ \ \ \left. -J_{\left( \frac{a+b}{2}\right) ^{+}}^{\alpha }\left( g\left(
f\circ h\right) \right) (b)-J_{\left( \frac{a+b}{2}\right) ^{-}}^{\alpha
}\left( g\left( f\circ h\right) \right) (a)\right.  \\ 
\\ 
\ \ \ \ -\frac{1}{\Gamma \left( \alpha \right) }\left.
\int\limits_{a}^{b}(h(b)-h(a))^{\alpha -1}g(t)f(h(t))h^{\prime
}(t)dt\right\}  \\ 
\\ 
=\frac{1}{\Gamma \left( \alpha \right) }\int\limits_{a}^{b}k(t)df(h(t))%
\end{array}%
\end{equation*}%
where%
\begin{equation*}
\begin{array}{c}
k(t)=\left\{ 
\begin{array}{ll}
\int\limits_{a}^{t}\left[ (h(b)-h(a))^{\alpha -1}-(h(b)-h(s))^{\alpha
-1}\right.  &  \\ 
\left. +(h(s)-h(a))^{\alpha -1}\right] g(s)h^{\prime }(s)ds & t\in \left[ a,%
\frac{a+b}{2}\right] , \\ 
\int\limits_{b}^{t}\left[ (h(b)-h(s))^{\alpha -1}-(h(s)-h(a))^{\alpha
-1}\right.  &  \\ 
\left. (h(b)-h(a))^{\alpha -1}\right] g(s)h^{\prime }(s)ds & t\in \left[ 
\frac{a+b}{2},b\right] .%
\end{array}%
\right. 
\end{array}%
\end{equation*}

\begin{proof}
It suffices to note that%
\begin{equation*}
\begin{array}{l}
I=\int\limits_{a}^{b}k(t)df(h(t)) \\ 
\\ 
=\int\limits_{a}^{\frac{a+b}{2}}\left\{ \int\limits_{a}^{t}\left[
(h(b)-h(a))^{\alpha -1}-(h(b)-h(s))^{\alpha -1}\right. \right.  \\ 
\\ 
\ \ \ \ +\left. \left. (h(s)-h(a))^{\alpha -1}\right] g(s)h^{\prime
}(s)ds\right\} df(h(t))%
\end{array}%
\end{equation*}%
\begin{equation*}
\begin{array}{l}
\ \ \ \ +\int\limits_{\frac{a+b}{2}}^{b}\left\{ \int\limits_{b}^{t}\left[
(h(b)-h(s))^{\alpha -1}-(h(s)-h(a))^{\alpha -1}\right. \right.  \\ 
\\ 
\ \ \ \ +\left. \left. (h(b)-h(a))^{\alpha -1}\right] g(s)h^{\prime
}(s)ds\right\} df(h(t)) \\ 
\\ 
=I_{1}+I_{2}.%
\end{array}%
\end{equation*}%
By integration by parts, we get%
\begin{equation*}
\begin{array}{l}
I_{1}=\int\limits_{a}^{\frac{a+b}{2}}\left\{ \int\limits_{a}^{t}\left[
(h(b)-h(a))^{\alpha -1}-(h(b)-h(s))^{\alpha -1}\right. \right.  \\ 
\\ 
\ \ \ \ +\left. \left. (h(s)-h(a))^{\alpha -1}\right] g(s)h^{\prime
}(s)ds\right\} df(h(t)) \\ 
\\ 
=\left. \left\{ \int\limits_{a}^{t}\left[ (h(b)-h(a))^{\alpha
-1}-(h(b)-h(s))^{\alpha -1}\right. \right. \right.  \\ 
\\ 
\ \ \ \ +\left. \left. \left. (h(s)-h(a))^{\alpha -1}\right] g(s)h^{\prime
}(s)ds\right\} f(h(t))\right\vert _{a}^{\frac{a+b}{2}}%
\end{array}%
\end{equation*}%
\begin{equation*}
\begin{array}{l}
\ \ \ \ -\left\{ \int\limits_{a}^{\frac{a+b}{2}}\left[ (h(b)-h(a))^{\alpha
-1}-(h(b)-h(t))^{\alpha -1}\right. \right.  \\ 
\\ 
\ \ \ \ +\left. \left. ((h(t)-h(a))^{\alpha -1}\right] g(t)f(h(t))h^{\prime
}(t)dt\right\}  \\ 
\\ 
=f\left( h\left( \frac{a+b}{2}\right) \right) \left\{ \int\limits_{a}^{\frac{%
a+b}{2}}\left[ (h(b)-h(a))^{\alpha -1}-(h(b)-h(s))^{\alpha -1}\right.
\right. 
\end{array}%
\end{equation*}%
\begin{equation*}
\begin{array}{l}
\ \ \ \ +\left. \left. (h(s)-h(a))^{\alpha -1}\right] g(s)h^{\prime
}(s)ds\right\}  \\ 
\\ 
\ \ \ \ -\left\{ \int\limits_{a}^{\frac{a+b}{2}}\left[ (h(b)-h(a))^{\alpha
-1}-(h(b)-h(t))^{\alpha -1}\right. \right.  \\ 
\\ 
\ \ \ \ +\left. \left. (h(t)-h(a))^{\alpha -1}\right] g(t)f(h(t))h^{\prime
}(t)dt\right\}  \\ 
\\ 
=\Gamma \left( \alpha \right) \left[ f\left( h\left( \frac{a+b}{2}\right)
\right) \left\{ J_{\left( \frac{a+b}{2}\right) ^{-}}^{\alpha
}g(h(a))-J_{a^{+}}^{\alpha }g(h(b))\right. \right. 
\end{array}%
\end{equation*}%
\begin{equation*}
\begin{array}{l}
\ \ \ \ +\left. \left. \frac{1}{\Gamma \left( \alpha \right) }%
\int\limits_{a}^{\frac{a+b}{2}}(h(b)-h(a))^{\alpha -1}g(s)h^{\prime
}(s)ds\right\} \right.  \\ 
\\ 
\ \ \ \ +\left. J_{a^{+}}^{\alpha }\left( g\left( f\circ h\right) \right)
(b)-J_{\left( \frac{a+b}{2}\right) ^{-}}^{\alpha }\left( g\left( f\circ
h\right) \right) (a)\right.  \\ 
\\ 
\ \ \ \ -\left. \frac{1}{\Gamma \left( \alpha \right) }\int\limits_{a}^{%
\frac{a+b}{2}}(h(b)-h(a))^{\alpha -1}g(t)f(h(t))h^{\prime }(t)dt\right] 
\end{array}%
\end{equation*}%
and similarly%
\begin{equation*}
\begin{array}{l}
I_{2}=\int\limits_{\frac{a+b}{2}}^{b}\left\{ \int\limits_{b}^{t}\left[
(h(b)-h(s))^{\alpha -1}-(h(s)-h(a))^{\alpha -1}\right. \right.  \\ 
\\ 
\ \ \ \ +\left. \left. (h(b)-h(a))^{\alpha -1}\right] g(s)h^{\prime
}(s)ds\right\} df(h(t)) \\ 
\\ 
=\left. \left\{ \int\limits_{b}^{t}\left[ (h(b)-h(s))^{\alpha
-1}-(h(s)-h(a))^{\alpha -1}\right. \right. \right.  \\ 
\\ 
\ \ \ \ +\left. \left. \left. (h(b)-h(a))^{\alpha -1}\right] g(s)h^{\prime
}(s)ds\right\} f(h(t))\right\vert _{\frac{a+b}{2}}^{b}%
\end{array}%
\end{equation*}%
\begin{equation*}
\begin{array}{l}
\ \ \ \ -\left\{ \int\limits_{\frac{a+b}{2}}^{b}\left[ (h(b)-h(t))^{\alpha
-1}-(h(t)-h(a))^{\alpha -1}\right. \right.  \\ 
\\ 
\ \ \ \ +\left. \left. (h(b)-h(a))^{\alpha -1}\right] g(t)f(h(t))h^{\prime
}(t)dt\right\}  \\ 
\\ 
=f\left( h\left( \frac{a+b}{2}\right) \right) \left\{ \int\limits_{\frac{a+b%
}{2}}^{b}\left[ (h(b)-h(s))^{\alpha -1}-(h(s)-h(a))^{\alpha -1}\right.
\right. 
\end{array}%
\end{equation*}%
\begin{equation*}
\begin{array}{l}
\ \ \ \ +\left. \left. (h(b)-h(a))^{\alpha -1}\right] g(s)h^{\prime
}(s)ds\right\}  \\ 
\\ 
\ \ \ \ -\left\{ \int\limits_{\frac{a+b}{2}}^{b}\left[ (h(b)-h(t))^{\alpha
-1}-(h(t)-h(a))^{\alpha -1}\right. \right.  \\ 
\\ 
\ \ \ \ +\left. \left. (h(b)-h(a))^{\alpha -1}\right] g(t)f(h(t))h^{\prime
}(t)dt\right\}  \\ 
\\ 
=\Gamma \left( \alpha \right) \left[ f\left( h\left( \frac{a+b}{2}\right)
\right) \left\{ J_{\left( \frac{a+b}{2}\right) ^{+}}^{\alpha
}g(h(b))-J_{b^{-}}^{\alpha }g(h(a))\right. \right. 
\end{array}%
\end{equation*}%
\begin{equation*}
\begin{array}{l}
\ \ \ \ +\left. \left. \frac{1}{\Gamma \left( \alpha \right) }\int\limits_{%
\frac{a+b}{2}}^{b}(h(b)-h(a))^{\alpha -1}g(s)h^{\prime }(s)ds\right\}
\right.  \\ 
\\ 
\ \ \ \ -\left\{ J_{\left( \frac{a+b}{2}\right) ^{+}}^{\alpha }\left(
g\left( f\circ h\right) \right) (b)+J_{b^{-}}^{\alpha }\left( g\left( f\circ
h\right) \right) (a)\right.  \\ 
\\ 
\ \ \ \ -\left. \left. \frac{1}{\Gamma \left( \alpha \right) }\int\limits_{%
\frac{a+b}{2}}^{b}(h(b)-h(a))^{\alpha -1}g(t)f(h(t))h^{\prime }(t)dt\right\} %
\right] .%
\end{array}%
\end{equation*}%
Thus, can write%
\begin{equation*}
\begin{array}{l}
I_{1}+I_{2}=\Gamma \left( \alpha \right) \left\{ f\left( h\left( \frac{a+b}{2%
}\right) \right) \left[ -J_{a^{+}}^{\alpha }g(h(b))-J_{b^{-}}^{\alpha
}g(h(a))\right. \right.  \\ 
\\ 
\ \ \ \ +\left. J_{\left( \frac{a+b}{2}\right) ^{+}}^{\alpha
}g(h(b))+J_{\left( \frac{a+b}{2}\right) ^{-}}^{\alpha }g(h(a))\right.  \\ 
\\ 
\ \ \ \ +\left. \left. \frac{1}{\Gamma \left( \alpha \right) }%
\int\limits_{a}^{\frac{a+b}{2}}(h(b)-h(a))^{\alpha -1}g(s)h^{\prime
}(s)ds\right. \right.  \\ 
\\ 
\ \ \ \ +\left. \frac{1}{\Gamma \left( \alpha \right) }\int\limits_{\frac{a+b%
}{2}}^{b}(h(b)-h(a))^{\alpha -1}g(s)h^{\prime }(s)ds\right] 
\end{array}%
\end{equation*}%
\begin{equation*}
\begin{array}{l}
\ \ \ \ +\left. J_{b^{-}}^{\alpha }\left( g\left( f\circ h\right) \right)
(a)+J_{a^{+}}^{\alpha }\left( g\left( f\circ h\right) \right) (b)\right.  \\ 
\\ 
\ \ \ \ \left. -J_{\left( \frac{a+b}{2}\right) ^{-}}^{\alpha }\left( g\left(
f\circ h\right) \right) (a)-J_{\left( \frac{a+b}{2}\right) ^{+}}^{\alpha
}\left( g\left( f\circ h\right) \right) (b)\right.  \\ 
\\ 
\ \ \ \ -\left. \frac{1}{\Gamma \left( \alpha \right) }\int\limits_{a}^{%
\frac{a+b}{2}}(h(b)-h(a))^{\alpha -1}g(t)f(h(t))h^{\prime }(t)dt\right.  \\ 
\\ 
\ \ \ \ -\left. \frac{1}{\Gamma \left( \alpha \right) }\int\limits_{\frac{a+b%
}{2}}^{b}(h(b)-h(a))^{\alpha -1}g(t)f(h(t))h^{\prime }(t)dt\right\} 
\end{array}%
\end{equation*}%
Multiplying the both sides $\left( \Gamma \left( \alpha \right) \right)
^{-1},\ $we obtain (\ref{2.4}) which complates the proof.
\end{proof}
\end{lemma}

\begin{remark}
\label{rem1} If we choose $h(x)=x\ $in Lemma \ref{lem2}, we have%
\begin{equation*}
\begin{array}{l}
\left\{ f\left( \frac{a+b}{2}\right) \left[ -J_{a^{+}}^{\alpha
}g(b)-J_{b^{-}}^{\alpha }g(a)+J_{\left( \frac{a+b}{2}\right) ^{+}}^{\alpha
}g(b)\right. \right. \\ 
\\ 
\ \ \ \ +\left. J_{\left( \frac{a+b}{2}\right) ^{-}}^{\alpha }g(a)+\dfrac{%
(b-a)^{\alpha -1}}{\Gamma \left( \alpha \right) }\int\limits_{a}^{b}g(s)ds%
\right] \\ 
\\ 
\ \ \ \ \left. +J_{a^{+}}^{\alpha }\left( gf\right) (b)+J_{b^{-}}^{\alpha
}\left( gf\right) (a)-J_{\left( \frac{a+b}{2}\right) ^{+}}^{\alpha }\left(
gf\right) (b)\right. \\ 
\\ 
\ \ \ \ \left. -J_{\left( \frac{a+b}{2}\right) ^{-}}^{\alpha }\left(
gf\right) (a)-\dfrac{(b-a)^{\alpha -1}}{\Gamma \left( \alpha \right) }%
\int\limits_{a}^{b}g(t)f(t)dt\right\} \\ 
\\ 
=\dfrac{1}{\Gamma \left( \alpha \right) }\int\limits_{a}^{b}k(t)df(t)%
\end{array}%
\end{equation*}%
where%
\begin{equation*}
k(t)=\left\{ 
\begin{array}{ll}
\int\limits_{a}^{t}\left[ (b-a))^{\alpha -1}-(b-s)^{\alpha -1}+(s-a)^{\alpha
-1}\right] g(s)ds & t\in \left[ a,\frac{a+b}{2}\right] , \\ 
\int\limits_{b}^{t}\left[ (b-s)^{\alpha -1}-(s-a)^{\alpha -1}+(b-a)^{\alpha
-1}\right] g(s)ds & t\in \left[ \frac{a+b}{2},b\right] .%
\end{array}%
\right.
\end{equation*}
\end{remark}

\begin{remark}
If we choose $h(x)=x\ $and$\ g(x)=1\ $in Lemma \ref{lem2}, we have%
\begin{equation*}
\begin{array}{l}
\left\{ f\left( \frac{a+b}{2}\right) \left[ (b-a)^{\alpha }\left( 1+\dfrac{1%
}{\alpha }\left( \dfrac{1}{2^{\alpha }}-2\right) \right) \right] \right. \\ 
\\ 
\ \ \ \ \left. +J_{a^{+}}^{\alpha }f(b)+J_{b^{-}}^{\alpha }f(a)-J_{\left( 
\frac{a+b}{2}\right) ^{+}}^{\alpha }f(b)\right. \\ 
\\ 
\ \ \ \ \left. -J_{\left( \frac{a+b}{2}\right) ^{-}}^{\alpha }f(a)-\dfrac{%
(b-a)^{\alpha -1}}{\Gamma \left( \alpha \right) }\int\limits_{a}^{b}f(t)dt%
\right\} \\ 
\\ 
=\dfrac{1}{\Gamma \left( \alpha \right) }\int\limits_{a}^{b}k(t)df(t)%
\end{array}%
\end{equation*}%
where%
\begin{equation*}
k(t)=\left\{ 
\begin{array}{ll}
\int\limits_{a}^{t}\left[ (b-a))^{\alpha -1}-(b-s)^{\alpha -1}+(s-a)^{\alpha
-1}\right] ds & t\in \left[ a,\frac{a+b}{2}\right] , \\ 
\int\limits_{b}^{t}\left[ (b-s)^{\alpha -1}-(s-a)^{\alpha -1}+(b-a)^{\alpha
-1}\right] ds & t\in \left[ \frac{a+b}{2},b\right] .%
\end{array}%
\right.
\end{equation*}
\end{remark}

\begin{remark}
If we choose $h(x)=x,\ \alpha =1\ $and$\ g(x)=1\ $in Lemma \ref{lem2}, we
obtain Lemma $2.1\ $in \cite{21}%
\begin{equation*}
\begin{array}{l}
\dfrac{1}{b-a}\int\limits_{a}^{b}f(t)dt-f\left( \dfrac{a+b}{2}\right) \\ 
=\left( b-a\right) \left( \int\limits_{0}^{\frac{1}{2}}tf^{\prime
}(ta+(1-t)b)dt+\int\limits_{\frac{1}{2}}^{1}\left( t-1\right) f^{\prime
}(ta+(1-t)b)dt\right) .%
\end{array}%
\end{equation*}

\begin{theorem}
\label{thm3} Let $f:I\rightarrow 
\mathbb{R}
$ be a differentiable mapping on $I^{\circ }$ and $f^{\prime }\in X_{h}^{p}%
\left[ a,b\right] $ with $a<b$ and $g:\left[ a,b\right] \rightarrow 
\mathbb{R}
$ is continuous. If $\left\vert f^{\prime }\right\vert \ $is convex on $%
\left[ a,b\right] ,\ $then the following inequality for fractional integrals
holds:%
\begin{equation*}
\begin{array}{l}
\left\{ f\left( \frac{a+b}{2}\right) \left[ -J_{a^{+}}^{\alpha
}g(b)-J_{b^{-}}^{\alpha }g(a)+J_{\left( \frac{a+b}{2}\right) ^{+}}^{\alpha
}g(b)\right. \right. \\ 
\\ 
\ \ \ \ \left. +J_{\left( \frac{a+b}{2}\right) ^{-}}^{\alpha }g(a)+\dfrac{%
(b-a)^{\alpha -1}}{\Gamma \left( \alpha \right) }\int\limits_{a}^{b}g(s)ds%
\right] \\ 
\\ 
\ \ \ \ \left. +J_{a^{+}}^{\alpha }\left( gf\right) (b)+J_{b^{-}}^{\alpha
}\left( gf\right) (a)-J_{\left( \frac{a+b}{2}\right) ^{+}}^{\alpha }\left(
gf\right) (b)\right. \\ 
\\ 
\ \ \ \ \left. -J_{\left( \frac{a+b}{2}\right) ^{-}}^{\alpha }\left(
gf\right) (a)-\dfrac{(b-a)^{\alpha -1}}{\Gamma \left( \alpha \right) }%
\int\limits_{a}^{b}g(t)f(t)dt\right\}%
\end{array}%
\end{equation*}%
\begin{equation*}
\begin{array}{l}
\leq \dfrac{\left\vert \left\vert g\right\vert \right\vert _{\left[ a,\frac{%
a+b}{2}\right] ,\infty }}{(b-a)\Gamma \left( \alpha \right) }\left\{
\int\limits_{a}^{\frac{a+b}{2}}\left[ \int\limits_{a}^{t}\left[
(b-a))^{\alpha -1}-(b-s)^{\alpha -1}+(s-a)^{\alpha -1}\right] ds\right]
\right. \\ 
\\ 
\ \ \ \ \times \left. \left( \left( b-t\right) \left\vert f^{\prime }\left(
a\right) \right\vert +\left( t-a\right) \left\vert f^{\prime }\left(
b\right) \right\vert \right) dt\right\} +\dfrac{\left\vert \left\vert
g\right\vert \right\vert _{\left[ \frac{a+b}{2},b\right] ,\infty }}{%
(b-a)\Gamma \left( \alpha \right) } \\ 
\\ 
\ \ \ \ \times \left\{ \int\limits_{\frac{a+b}{2}}^{b}\left[
\int\limits_{b}^{t}\left[ (b-s)^{\alpha -1}-(s-a)^{\alpha -1}+(b-a)^{\alpha
-1}\right] g(s)ds\right] \right. \\ 
\\ 
\ \ \ \ \times \left. \left( \left( b-t\right) \left\vert f^{\prime }\left(
a\right) \right\vert +\left( t-a\right) \left\vert f^{\prime }\left(
b\right) \right\vert \right) dt\right\}%
\end{array}%
\end{equation*}%
with $\alpha >0.$

\begin{proof}
If $\left\vert f^{\prime }\right\vert \ $is convex on $\left[ a,b\right] ,\ $%
we know that for $t\in \left[ a,b\right] $%
\begin{equation*}
\begin{array}{l}
\left\vert f^{\prime }(t)\right\vert = \\ 
\\ 
\left\vert f^{\prime }\left( \dfrac{b-t}{b-a}a+\dfrac{t-a}{b-a}b\right)
\right\vert \leq \dfrac{b-t}{b-a}\left\vert f^{\prime }\left( a\right)
\right\vert +\dfrac{t-a}{b-a}\left\vert f^{\prime }\left( b\right)
\right\vert .%
\end{array}%
\end{equation*}%
From Remark \ref{rem1}, we have%
\begin{equation*}
\begin{array}{l}
\left\{ f\left( \frac{a+b}{2}\right) \left[ -J_{a^{+}}^{\alpha
}g(b)-J_{b^{-}}^{\alpha }g(a)+J_{\left( \frac{a+b}{2}\right) ^{+}}^{\alpha
}g(b)\right. \right.  \\ 
\\ 
\ \ \ \ \left. +J_{\left( \frac{a+b}{2}\right) ^{-}}^{\alpha }g(a)+\dfrac{%
(b-a)^{\alpha -1}}{\Gamma \left( \alpha \right) }\int\limits_{a}^{b}g(s)ds%
\right] 
\end{array}%
\end{equation*}%
\begin{equation*}
\begin{array}{l}
\ \ \ \ +J_{a^{+}}^{\alpha }\left( gf\right) (b)+J_{b^{-}}^{\alpha }\left(
gf\right) (a)-J_{\left( \frac{a+b}{2}\right) ^{+}}^{\alpha }\left( gf\right)
(b) \\ 
\\ 
\ \ \ \ \left. -J_{\left( \frac{a+b}{2}\right) ^{-}}^{\alpha }\left(
gf\right) (a)-\dfrac{(b-a)^{\alpha -1}}{\Gamma \left( \alpha \right) }%
\int\limits_{a}^{b}g(t)f(t)dt\right\} 
\end{array}%
\end{equation*}%
\begin{equation*}
\begin{array}{l}
\leq \frac{1}{\Gamma \left( \alpha \right) }\left\{ \int\limits_{a}^{\frac{%
a+b}{2}}\left\vert \int\limits_{a}^{t}\left[ (b-a))^{\alpha
-1}-(b-s)^{\alpha -1}+(s-a)^{\alpha -1}\right] g(s)ds\right\vert \left\vert
f^{\prime }\left( t\right) \right\vert dt\right.  \\ 
\ \ \ \ \left. +\int\limits_{\frac{a+b}{2}}^{b}\left\vert \int\limits_{b}^{t}%
\left[ (b-s)^{\alpha -1}-(s-a)^{\alpha -1}+(b-a)^{\alpha -1}\right]
ds\right\vert g(s)\left\vert f^{\prime }\left( t\right) \right\vert
dt\right\}  \\ 
\leq \dfrac{1}{(b-a)\Gamma \left( \alpha \right) }\left\{ \int\limits_{a}^{%
\frac{a+b}{2}}\left[ \int\limits_{a}^{t}\left[ (b-a))^{\alpha
-1}-(b-s)^{\alpha -1}+(s-a)^{\alpha -1}\right] g(s)ds\right] \right. 
\end{array}%
\end{equation*}%
\begin{equation*}
\begin{array}{l}
\ \ \ \ \times \left( \left( b-t\right) \left\vert f^{\prime }\left(
a\right) \right\vert +\left( t-a\right) \left\vert f^{\prime }\left(
b\right) \right\vert \right) dt \\ 
\ \ \ \ +\int\limits_{\frac{a+b}{2}}^{b}\left[ \int\limits_{b}^{t}\left[
(b-s)^{\alpha -1}-(s-a)^{\alpha -1}+(b-a)^{\alpha -1}\right] g(s)ds\right] 
\end{array}%
\end{equation*}%
\begin{equation*}
\begin{array}{l}
\ \ \ \ \times \left. \left( \left( b-t\right) \left\vert f^{\prime }\left(
a\right) \right\vert +\left( t-a\right) \left\vert f^{\prime }\left(
b\right) \right\vert \right) dt\right\}  \\ 
\leq \dfrac{\left\vert \left\vert g\right\vert \right\vert _{\left[ a,\frac{%
a+b}{2}\right] ,\infty }}{(b-a)\Gamma \left( \alpha \right) }\left\{
\int\limits_{a}^{\frac{a+b}{2}}\left[ \int\limits_{a}^{t}\left[
(b-a))^{\alpha -1}-(b-s)^{\alpha -1}+(s-a)^{\alpha -1}\right] ds\right]
\right.  \\ 
\ \ \ \ \times \left. \left( \left( b-t\right) \left\vert f^{\prime }\left(
a\right) \right\vert +\left( t-a\right) \left\vert f^{\prime }\left(
b\right) \right\vert \right) dt\right\}  \\ 
\ \ \ \ +\dfrac{\left\vert \left\vert g\right\vert \right\vert _{\left[ 
\frac{a+b}{2},b\right] ,\infty }}{(b-a)\Gamma \left( \alpha \right) }\left\{
\int\limits_{\frac{a+b}{2}}^{b}\left[ \int\limits_{b}^{t}\left[
(b-s)^{\alpha -1}-(s-a)^{\alpha -1}+(b-a)^{\alpha -1}\right] g(s)ds\right]
\right.  \\ 
\ \ \ \ \times \left. \left( \left( b-t\right) \left\vert f^{\prime }\left(
a\right) \right\vert +\left( t-a\right) \left\vert f^{\prime }\left(
b\right) \right\vert \right) dt\right\} 
\end{array}%
\end{equation*}%
\begin{equation*}
\begin{array}{l}
\leq \dfrac{\left\vert \left\vert g\right\vert \right\vert _{\left[ a,\frac{%
a+b}{2}\right] ,\infty }}{(b-a)\Gamma \left( \alpha \right) }\left\{
\int\limits_{a}^{\frac{a+b}{2}}\left[ \int\limits_{a}^{t}\left[
(b-a))^{\alpha -1}-(b-s)^{\alpha -1}+(s-a)^{\alpha -1}\right] ds\right]
\right.  \\ 
\ \ \ \ \times \left. \left( \left( b-t\right) \left\vert f^{\prime }\left(
a\right) \right\vert +\left( t-a\right) \left\vert f^{\prime }\left(
b\right) \right\vert \right) dt\right\}  \\ 
\ \ \ \ +\dfrac{\left\vert \left\vert g\right\vert \right\vert _{\left[ 
\frac{a+b}{2},b\right] ,\infty }}{(b-a)\Gamma \left( \alpha \right) }\left\{
\int\limits_{\frac{a+b}{2}}^{b}\left[ \int\limits_{b}^{t}\left[
(b-s)^{\alpha -1}-(s-a)^{\alpha -1}+(b-a)^{\alpha -1}\right] g(s)ds\right]
\right.  \\ 
\ \ \ \ \times \left. \left( \left( b-t\right) \left\vert f^{\prime }\left(
a\right) \right\vert +\left( t-a\right) \left\vert f^{\prime }\left(
b\right) \right\vert \right) dt\right\} .%
\end{array}%
\end{equation*}%
This completes the proof.
\end{proof}
\end{theorem}
\end{remark}


\begin{thebibliography}{99}
\bibitem{1} A.A. Kilbas, H.M. Srivastava, J.J. Trujillo, \textit{Theory and
Applications of Fractional Differantial Equations}, Elsevier B.V.,
Amsterdam, 2006.\bigskip

\bibitem{2} A. Akkurt, H. Y\i ld\i r\i m, \textit{On Feng Qi Type Integral
Inequalities For Generalized Fractional Integrals}, IAAOJ, Scientific
Science, \textbf{1(2)}, (2013), 22-25.\bigskip

\bibitem{3} A. Akkurt, Z. K\i rtay, H. Y\i ld\i r\i m, \textit{Generalized
Fractional Integrals Inequalities for Continuous Random Variables}, Journal
of Probability and Statistics, Volume \textbf{2015} (2015), Article ID
958980, 7 pages, http://dx.doi.org/10.1155/2015/958980.\bigskip

\bibitem{4} E. F. Beckenbach, \textit{Convex functions}, Bull. Amer. Math.
Soc., \textbf{54} (1948) 439-460.
http://dx.doi.org/10.1090/s0002-9904-1948-08994-7\bigskip

\bibitem{5} E. Set, \.{I}. \.{I}\c{s}can, M. E. \"{O}zdemir and M. Z. Sar\i
kaya, \textit{On new Hermite-Hadamard-Fejer type inequalities for convex
functions via fractional integrals}, Applied Mathematics and Computation, 
\textbf{259}(2015), 875--881.\bigskip

\bibitem{6} H. Y\i ld\i r\i m, Z. K\i rtay, \textit{Ostrowski Inequality for
Generalized Fractional Integral and Related Inequalities}, Malaya Journal of
Matematik: Volume \textbf{2}, Issue 3, 2014, pp. 322-329.\bigskip

\bibitem{7} \.{I}. \.{I}\c{s}can, \textit{Hermite--Hadamard--Fej\'{e}r type
inequalities for convex functions via fractional integrals}, Stud. Univ. Babe%
\c{s}-Bolyai Math. \textbf{60 }(2015), No. 3, 355-366\bigskip

\bibitem{8} J. Park, \textit{Inequalities of Hermite-Hadamard-Fej\'{e}r Type
for Convex Functions via Fractional Integrals},International Journal of
Mathematical Analysis,Vol. \textbf{8}, 2014, no. 59, pp. 2927-2937,
http://dx.doi.org/10.12988/ijma.2014.411378.\bigskip

\bibitem{9} K. L. Tseng, S. R. Hwang, S. S. Dragomir, \textit{Fej\'{e}r-type
inequalities (I)}, J. Inequal. Appl., \textbf{2010} (2010), Art ID: 531976,
7 pages. http://dx.doi.org/10.1155/2010/531976.\bigskip

\bibitem{10} L. Fej\'{e}r, \textit{Uberdie Fourierreihen, II}, Math.
Naturwise. Anz Ungar. Akad., Wiss \textbf{24} (1906) 369--390 (in
Hungarian).\bigskip

\bibitem{11} M.Z. Sarikaya, E. Set, H. Yald\i z, N. Ba\c{s}ak, \textit{%
Hermite--Hadamard's inequalities for fractional integrals and related
fractional inequalities}, Math. Comput.Modell. \textbf{57} (9) (2013)
2403--2407.\bigskip

\bibitem{12} M. E. \"{O}zdemir, M. Avic, H. Kavurmaci, \textit{%
Hermite-Hadamard type inequalities for s-convex and s-concave functions via
fractional integrals}, Turkish Journal of Science \textbf{1} (2016):
28-40.\bigskip

\bibitem{13} M. Z. Sarikaya, \textit{On new Hermite Hadamard Fej\'{e}r type
integral inequalities}, Stud. Univ. Babes-Bolyai Math. \textbf{57(3)} (2012)
377-386.\bigskip

\bibitem{14} M. Tun\c{c}, \textit{On some new inequalities for convex
functions}, Turk. J. Math., \textbf{35} (2011), 1-7.\bigskip

\bibitem{15} M. Z. Sarikaya, S. Erden, \textit{On the Hermite-Hadamard-Fej%
\'{e}r type integral inequality for convex function}, Turk. J. Anal. Number
Theory \textbf{2(3)} (2014) 85--89.\bigskip

\bibitem{16} M. Z. Sarikaya, A. Sa\u{g}lam and H. Y\i ld\i r\i m, \textit{On
some Hadamard-type inequalities for h-convex functions}, J. Math. Ineq., 
\textbf{2(3) }(2008), 335-341.\bigskip

\bibitem{17} S.G. Samko, A.A. Kilbas, and O.I. Marichev, \textit{Fractional
Integrals and Derivatives - Theory and Applications}, Gordon and Breach,
Linghorne, 1993.\bigskip

\bibitem{18} S. Belarbi, Z. Dahmani, \textit{On some new fractional integral
inequalities}, J. Ineq. Pure Appl. Math. \textbf{10} (3) (2009). Art.
86.\bigskip

\bibitem{19} S. K\i l\i n\c{c}, H. Y\i ld\i r\i m, \textit{Generalized
Fractional Integral Inequalities Involving Hipergeometric Operators},
IJPAM., vol.\textbf{101} no.1 (2015) pp.71-82.\bigskip

\bibitem{20} Z. Dahmani, \textit{On Minkowski and Hermite-Hadamard integral
inequalities via fractional integration}, Ann. Funct. Anal. \textbf{1(1)}
(2010) 51-58.\bigskip

\bibitem{21} U.S. K\i rmac\i , \textit{Inequalities for differentiable
mappings and applications to special means of real numbers and to midpoint
formula}, Appl. Math. Comput. \textbf{147} (1) (2004) 137--146.
\end{thebibliography}
\end{document}